# On the approximation of mean densities of random closed sets

LUIGI AMBROSIO[1], VINCENZO CAPASSO[2,3] and ELENA VILLA[2]

[1]*Scuola Normale Superiore, p.za dei Cavalieri 7, 56126 Pisa, Italy.* E-mail: *l.ambrosio@sns.it*
[2]*Department of Mathematics, University of Milan, via Saldini 50, 20133 Milano, Italy.*
E-mail: *vincenzo.capasso@unimi.it*
[3]*ADAMSS (Centre for Advanced Applied Mathematical and Statistical Sciences), University of Milan, Italy.* E-mail: *elena.villa@unimi.it*

Many real phenomena may be modelled as random closed sets in $\mathbb{R}^d$, of different Hausdorff dimensions. In many real applications, such as fiber processes and $n$-facets of random tessellations of dimension $n \leq d$ in spaces of dimension $d \geq 1$, several problems are related to the estimation of such mean densities. In order to confront such problems in the general setting of spatially inhomogeneous processes, we suggest and analyze an approximation of mean densities for sufficiently regular random closed sets. We show how some known results in literature follow as particular cases. A series of examples throughout the paper are provided to illustrate various relevant situations.

*Keywords:* mean densities; random measures; stochastic geometry

## 1. Introduction

Many real phenomena may be modelled as random closed sets in $\mathbb{R}^d$, of different Hausdorff dimensions (see, e.g., [3, 4, 13, 17, 19, 23, 24]).

We recall that a *random closed set* $\Theta$ in $\mathbb{R}^d$ is a measurable map

$$\Theta : (\Omega, \mathcal{F}, \mathbb{P}) \longrightarrow (\mathbb{F}, \sigma_{\mathbb{F}}),$$

where $\mathbb{F}$ denotes the class of closed subsets in $\mathbb{R}^d$ and $\sigma_{\mathbb{F}}$ is the $\sigma$-algebra generated by the so-called *hit-or-miss topology* (see [16]).

Let $\Theta_n$ be almost surely (a.s.) a set of locally finite Hausdorff $n$-dimensional measure and denote by $\mathcal{H}^n$ the $n$-dimensional Hausdorff measure on $\mathbb{R}^d$. The set $\Theta_n$ induces a random measure $\mu_{\Theta_n}$ defined by

$$\mu_{\Theta_n}(A) := \mathcal{H}^n(\Theta_n \cap A), \qquad A \in \mathcal{B}_{\mathbb{R}^d}$$

(for a discussion of the delicate issue of measurability of the random variables $\mathcal{H}^n(\Theta_n \cap A)$, we refer to [2, 26]).

Under suitable regularity assumptions on a random closed set $\Theta_n \subset \mathbb{R}^d$ with locally finite $n$-dimensional measure a.s., in [7] the concept of *mean density of $\Theta_n$*, that is, the density of the expected measure

$$\mathbb{E}[\mu_{\Theta_n}](A) := \mathbb{E}[\mathcal{H}^n(\Theta_n \cap A)], \qquad A \in \mathcal{B}_{\mathbb{R}^d}, \tag{1}$$

with respect to the standard Lebesgue measure $\nu^d$ on $\mathbb{R}^d$, has been revisited in terms of expected values of a suitable class of linear functionals (Delta functions *à* la Dirac).

It is clear that if $n < d$ and $\mu_{\Theta_n(\omega)}$ is a Radon measure for almost every $\omega \in \Omega$, then it is singular with respect to the $d$-dimensional Lebesgue measure $\nu^d$. On the other hand, in dependence of the probability law of $\Theta_n$, the expected measure may be either singular or absolutely continuous with respect to $\nu^d$; in the latter case, we say that a random closed set is *absolutely continuous in mean*, and its mean density is the classical Radon–Nikodym derivative of $\mathbb{E}[\mu_{\Theta_n}]$ with respect to $\nu^d$.

In many real applications, such as fiber processes and $n$-facets of random tessellations of dimension $n \leq d$ in spaces of dimension $d \geq 1$, several problems are related to the estimation of their mean densities (see, e.g., [3, 5, 12, 18, 19, 24]). The scope of this paper is to provide concrete approximations that may lead to (possibly local) estimators of mean densities of lower-dimensional inhomogeneous random closed sets; the globally stationary case has been the subject of various papers and books in the literature (see, e.g., [3]), while the problem of local density estimation for $n = 0$ (absolutely continuous random variables) has been largely solved for quite some time – this now being part of the standard literature – by means of either histograms, or kernel estimators (see, e.g., [10, 22]). In order to confront such problems in the general setting of spatially inhomogeneous processes, we suggest and analyze here an approximation of mean densities for sufficiently regular random closed sets $\Theta_n$ in $\mathbb{R}^d$. We shall show how some known results in the literature follow as particular cases.

A series of examples throughout the paper are provided to illustrate various relevant situations.

## 2. Preliminaries and notation

In this section, we collect some basic facts and terminology that will be useful in the sequel.

We shall define a *Radon measure* in $\mathbb{R}^d$ to be any non-negative and $\sigma$-additive set function $\mu$ defined on the Borel $\sigma$-algebra $\mathcal{B}_{\mathbb{R}^d}$ which is finite on bounded sets.

We know that every Radon measure $\mu$ on $\mathbb{R}^d$ can be represented in the form

$$\mu = \mu_{\ll} + \mu_{\perp},$$

where $\mu_{\ll}$ and $\mu_{\perp}$ are the absolutely continuous part of $\mu$ with respect to $\nu^d$ and the singular part of $\mu$, respectively. Denoting by $B_r(x)$ the closed ball with center $x$ and radius $r$, as a consequence of the Besicovitch derivation theorem (see [1], page 54), we have that the limit

$$\delta\mu(x) := \lim_{r\downarrow 0} \frac{\mu(B_r(x))}{\nu^d(B_r(x))}$$

exists in $\mathbb{R}$ for $\nu^d$-a.e. $x \in \mathbb{R}^d$ and it is a version of the Radon–Nikodym derivative of $\mu_\ll$, while $\mu_\perp$ is the restriction of $\mu$ to the $\nu^d$-negligible set $\{x \in \mathbb{R}^d : \delta\mu(x) = \infty\}$. We call $\delta\mu$ the *density* of the measure $\mu$ with respect to $\nu^d$.

Let us consider a random closed set in $\mathbb{R}^d$ with locally finite $n$-dimensional measure a.s., say $\Theta_n$, such that the expected measure $\mathbb{E}[\mu_{\Theta_n}]$ induced by $\Theta_n$ defined by (1) is a Radon measure and let $\delta\mathbb{E}[\mu_{\Theta_n}]$ be its density. We refer to $\delta\mathbb{E}[\mu_{\Theta_n}]$ as the *mean density of* $\Theta_n$. (For a more detailed discussion about the density of deterministic and random closed sets, see [7].) For any random closed set $\Theta_n$ in $\mathbb{R}^d$ with locally finite $n$-dimensional measure for some $n < d$ a.s., while it is clear that $\mu_{\Theta_n(\omega)}$ is a singular measure, it may well happen (e.g., when the process $\Theta_n$ is stationary) that the expected measure $\mathbb{E}[\mu_{\Theta_n}]$ is absolutely continuous with respect to $\nu^d$, so that $\delta\mathbb{E}[\mu_{\Theta_n}]$ is the Radon–Nikodym density of $\mathbb{E}[\mu_{\Theta_n}]$. In this case, the following definition is given [6].

***Definition 1 (Absolute continuity in mean).*** *Let $\Theta_n$ be a random closed set in $\mathbb{R}^d$ with locally finite $n$-dimensional measure a.s. such that its associated expected measure $\mathbb{E}[\mu_{\Theta_n}]$ is a Radon measure. We say that $\Theta_n$ is absolutely continuous in mean if $\mathbb{E}[\mu_{\Theta_n}]$ is absolutely continuous with respect to $\nu^d$. In this case, we denote by $\lambda_{\Theta_n}$ its density.*

***Remark 2 (The $0$-dimensional and $d$-dimensional cases).*** If $n = 0$ and $\Theta_0(\omega) = \{X(\omega)\}$ is a random point, then $\mathbb{E}[\mathcal{H}^0(\Theta_0 \cap \cdot)] = \mathbb{P}(X \in \cdot)$. Therefore, $\Theta_0$ is absolutely continuous in mean if and only if the law of $X$ is absolutely continuous and $\lambda_{\Theta_0}$ coincides with the probability density function of $X$.

On the other hand, $\mathbb{E}[\mu_{\Theta_d}]$ is always absolutely continuous with respect to $\nu^d$ and by a simple application of Fubini's theorem (in $\Omega \times \mathbb{R}^d$, with the product measure $\mathbb{P} \times \nu^d$), we have that [21] (see also [15])

$$\mathbb{E}[\mu_{\Theta_d}(B)] = \int_B \mathbb{P}(x \in \Theta_d)\,\mathrm{d}x \qquad \forall B \in \mathcal{B}_{\mathbb{R}^d},$$

so $\lambda_{\Theta_d}(x) = \mathbb{P}(x \in \Theta_d)$ for $\nu^d$-a.e. $x \in \mathbb{R}^d$.

We mention that a stronger definition of absolute continuity in mean of a random closed set has been given in [6] in terms of the expected measure of its boundary.

The aim of the present paper is to provide an approximation of the expected measure $\mathbb{E}[\mu_{\Theta_n}]$ under quite general regularity assumptions on the random closed set $\Theta_n$, which provides, as a by-product, an approximation of the mean density $\lambda_{\Theta_n}$ when $\mathbb{E}[\mu_{\Theta_n}]$ is absolutely continuous with respect to $\nu^d$. To this end, the notion of *Minkowski content* of closed sets in $\mathbb{R}^d$ plays here a fundamental role. Denoting by $b_k$ the volume of the unit ball in $\mathbb{R}^k$ and by $S_{\oplus r}$ the closed $r$-neighborhood of $S$, that is,

$$S_{\oplus r} := \{x \in \mathbb{R}^d : \exists y \in S \text{ with } |x - y| \le r\}$$

(known also as the *Minkowski enlargement* of $S$ with the closed ball $B_r(0)$), the $n$-dimensional *Minkowski content* of a closed set $S \subset \mathbb{R}^d$ is defined by

$$\lim_{r \downarrow 0} \frac{\nu^d(S_{\oplus r})}{b_{d-n} r^{d-n}},$$

whenever the limit exists.

General results concerning the existence of the Minkowski content of closed subsets in $\mathbb{R}^d$ are well known in the literature, related to rectifiability properties of the sets involved. Given an integer $0 \leq n \leq d$, we say that a set $C \subset \mathbb{R}^d$ is *countably $\mathcal{H}^n$-rectifiable* if there exist countably many Lipschitz maps $f_i \colon \mathbb{R}^n \to \mathbb{R}^d$ such that

$$\mathcal{H}^n\left(C \setminus \bigcup_{i=1}^{\infty} f_i(\mathbb{R}^n)\right) = 0.$$

Rectifiable sets include piecewise $C^1$ sets, and they still have nice properties from the measure-theoretic viewpoint (e.g., one can define an $n$-dimensional tangent space to them, in an approximate sense); we refer to [1] for the basic properties of this class of sets. We quote the following result from [1], page 110.

**Theorem 3.** *Let $S \subset \mathbb{R}^d$ be a countably $\mathcal{H}^n$-rectifiable compact set and assume that*

$$\eta(B_r(x)) \geq \gamma r^n \qquad \forall x \in S, \forall r \in (0, 1) \tag{2}$$

*holds for some $\gamma > 0$ and some Radon measure $\eta$ in $\mathbb{R}^d$ absolutely continuous with respect to $\mathcal{H}^n$. Then*

$$\lim_{r \downarrow 0} \frac{\nu^d(S_{\oplus r})}{b_{d-n} r^{d-n}} = \mathcal{H}^n(S).$$

In particular, in Lemma 6, we shall prove a local version of the above theorem.

## 3. Approximation of mean densities

Throughout the paper, we consider countably $\mathcal{H}^n$-rectifiable random closed sets in $\mathbb{R}^d$ (i.e., $\mathbb{P}$-a.e. $\omega \in \Omega$, $\Theta_n(\omega) \subset \mathbb{R}^d$ is a countably $\mathcal{H}^n$-rectifiable closed set), with Radon measure $\mathbb{E}[\mu_{\Theta_n}]$.

By the aforementioned Besicovitch derivation theorem, a natural approximation at the scale $r$ of $\lambda_{\Theta_n}$ can be given by

$$\delta^{(r)}\mathbb{E}[\mu_{\Theta_n}](x) := \frac{\mathbb{E}[\mu_{\Theta_n}](B_r(x))}{b_d r^d} = \mathbb{E}\left[\frac{\mu_{\Theta_n}(B_r(x))}{b_d r^d}\right].$$

It can be easily proven that the associated measures $\delta^{(r)}\mathbb{E}[\mu_{\Theta_n}]\nu^d$ weakly-* converge to the measure $\mathbb{E}[\mu_{\Theta_n}]$ as $r \downarrow 0$ (i.e., in the duality with continuous and compactly supported functions). This convergence result can also be understood by noting that $\delta^{(r)}\mathbb{E}[\mu_{\Theta_n}](x)$ is the convolution of the measure $\mathbb{E}[\mu_{\Theta_n}]$ with the kernels (here, $\mathbf{1}_E$ stands for the characteristic function of $E$)

$$\rho_r(y) := \frac{1}{b_d r^d}\mathbf{1}_{B_r(0)}(y).$$

Noting that a point $y$ belongs to $B_r(x)$ if and only if the point $x$ belongs to the Minkowski enlargement $y_{\oplus r}$ of $y$, it follows that whenever $\Theta_0$ is a random point in $\mathbb{R}^d$,

$$\delta^{(r)}\mathbb{E}[\mu_{\Theta_0}](x) = \frac{\mathbb{E}[\mathcal{H}^0(\Theta_0 \cap B_r(x))]}{b_d r^d} = \frac{\mathbb{E}[\mathbf{1}_{\Theta_{0\oplus r}}(x)]}{b_d r^d}.$$

In the particular case $d=1$ with $\Theta_0 = X$ a random variable, we have that

$$\delta^{(r)}\mathbb{E}[\mu_X](x) = \frac{\mathbb{P}(X \in [x-r, x+r])}{2r}; \tag{3}$$

if $X$ is a random variable with absolutely continuous law and pdf $f_X$, we have (see Remark 2) that $\lambda_X = f_X$, so that (3) leads to the usual histogram estimation of probability densities (see [20], Section VII.13).

We point out that even if a natural approximation of $\mathbb{E}[\mu_{\Theta_n}]$ can be given by the functions $\delta^{(r)}\mathbb{E}[\mu_{\Theta_n}]$ defined above, problems might arise in the estimation of $\mathbb{E}[\mu_{\Theta_n}(B_r(x))]$, as the computation of the Hausdorff measure $\mathcal{H}^n(\Theta_n(\omega) \cap B_r(x))$ is typically non-trivial. As a matter of fact, a computer graphics representation of lower-dimensional sets in $\mathbb{R}^2$ is anyway provided in terms of pixels, which can only offer a 2-D box approximation of points in $\mathbb{R}^2$ (an interesting discussion on this is contained in [14]). Therefore, we are led to consider a new approximation, based on the Lebesgue measure (much more robust and computable) of the Minkowski enlargement of the random set in question. This procedure is obviously consistent with (3). We now state the main theorem of this paper.

Given a random closed set $\Theta_n$ with locally finite $n$-dimensional measure a.s. in $\mathbb{R}^d$ (the particular case $n=d$ is trivial), for any compact window $W \subset \mathbb{R}^d$, let $\Gamma_W(\Theta_n): \Omega \longrightarrow \mathbb{R}$ be the function defined as follows:

$$\Gamma_W(\Theta_n) := \sup\{\gamma \geq 0 : \exists \text{ a probability measure } \eta \ll \mathcal{H}^n \text{ such that}$$
$$\eta(B_r(x)) \geq \gamma r^n \ \forall x \in \Theta_n \cap W_{\oplus 1}, \forall r \in (0,1)\}. \tag{4}$$

**Theorem 4 (Main result).** *Let $\Theta_n$ be a countably $\mathcal{H}^n$-rectifiable random closed set in $\mathbb{R}^d$ such that $\mathbb{E}[\mu_{\Theta_n}]$ is a Radon measure. Assume that for any compact window $W \subset \mathbb{R}^d$, there exists a random variable $Y$ with $\mathbb{E}[Y] < \infty$, such that $1/\Gamma_W(\Theta_n) \leq Y$ a.s. Then*

$$\lim_{r \downarrow 0} \int_A \frac{\mathbb{P}(x \in \Theta_{n \oplus r})}{b_{d-n} r^{d-n}}\, \mathrm{d}x = \mathbb{E}[\mu_{\Theta_n}](A) \tag{5}$$

*for any bounded Borel set $A \subset \mathbb{R}^d$ such that*

$$\mathbb{E}[\mu_{\Theta_n}](\partial A) = 0. \tag{6}$$

In particular, if $\Theta_n$ is absolutely continuous in mean, we have

$$\lim_{r\downarrow 0}\int_A \frac{\mathbb{P}(x\in\Theta_{n\oplus r})}{b_{d-n}r^{d-n}}\,\mathrm{d}x = \int_A \lambda_{\Theta_n}(x)\,\mathrm{d}x \qquad (7)$$

for any bounded Borel set $A\subset\mathbb{R}^d$ with $\nu^d(\partial A)=0$, where $\lambda_{\Theta_n}$ is the mean density of $\Theta_n$. Finally, if $\Theta_n$ is stationary, we have

$$\lim_{r\downarrow 0}\frac{\mathbb{P}(x_0\in\Theta_{n\oplus r})}{b_{d-n}r^{d-n}} = L_{\Theta_n}\in\mathbb{R}_+ \qquad \forall x_0\in\mathbb{R}^d, \qquad (8)$$

having denoted by $L_{\Theta_n}$ the constant value of $\lambda_{\Theta_n}$.

***Remark 5.*** In the statement of the above theorem, we introduced the auxiliary random variable $Y$ in order to avoid the non-trivial issue of the measurability of $\Gamma_W$; as a matter of fact, in all examples, one can estimate $1/\Gamma_W$ from above in a measurable way.

Note also that if $\Theta_n$ satisfies the assumption of the theorem for some closed $W$, then it satisfies the assumption for all closed $W'\subset W$; analogously, any random closed set $\Theta'_n$ contained almost surely in $\Theta_n$ still satisfies the assumption of the theorem.

Finally, note that condition (6), when restricted to bounded open sets $A$, is "generically satisfied" in the following sense: given any family of bounded open sets $\{A_t\}_{t\in\mathbb{R}}$ with $\mathrm{clos}\,A_s\subseteq A_t$ for $s<t$ ($\mathrm{clos}\,A$ denotes here the closure of the set $A$), the set

$$T:=\{t\in\mathbb{R}:\mathbb{E}[\mu_{\Theta_n}](\partial A_t)>0\}$$

is at most countable. This is due to the fact that the sets $\{\partial A_t\}_{t\in T}$ are pairwise disjoint and all with strictly positive $\mathbb{E}[\mu_{\Theta_n}]$-measure (in particular, the sets $\{t\in(-\infty,m):\mathbb{E}[\mu_{\Theta_n}](\partial A_t)>1/m\}$ have cardinality at most $m\mathbb{E}[\mu_{\Theta_n}](A_m)$).

### 3.1. Proof of the main result

Note that by a straightforward application of Fubini's theorem, the limit in (5) is equivalent to

$$\lim_{r\downarrow 0}\frac{\mathbb{E}[\nu^d(\Theta_{n\oplus r}\cap A)]}{b_{d-n}r^{d-n}} = \mathbb{E}[\mathcal{H}^n(\Theta_n\cap A)], \qquad (9)$$

which might be regarded as the *local mean n-dimensional Minkowski content* of $\Theta_n$. The problem hence reduces to finding conditions on $\Theta_n$ which ensure that (9) holds. To this end, let us observe that if $\Theta_n$ is such that almost every realization $\Theta_n(\omega)$ has $n$-dimensional Minkowski content equal to the $n$-dimensional Hausdorff measure, that is,

$$\lim_{r\downarrow 0}\frac{\nu^d(\Theta_{n\oplus r}(\omega))}{b_{d-n}r^{d-n}} = \mathcal{H}^n(\Theta_n(\omega)), \qquad (10)$$

then it is clear that, taking the expected values on both sides, (9) is strictly related to the possibility of exchanging limit and expectation. We therefore ask whether (10) implies a similar result when we consider the intersection of $\Theta_{n_{\oplus r}}(\omega)$ with an open set $A$ in $\mathbb{R}^d$, and for which kind of random closed sets the convergence above is dominated, so that exchanging limit and expectation is allowed.

The following result is a local version of Theorem 3.

**Lemma 6.** *Let $S$ be a compact subset of $\mathbb{R}^d$ satisfying the hypotheses of Theorem 3. Then, for any $A \in \mathcal{B}_{\mathbb{R}^d}$ such that*

$$\mathcal{H}^n(S \cap \partial A) = 0, \tag{11}$$

*the following holds:*

$$\lim_{r \downarrow 0} \frac{\nu^d(S_{\oplus r} \cap A)}{b_{d-n} r^{d-n}} = \mathcal{H}^n(S \cap A). \tag{12}$$

**Proof.** If $n = d$, then equality (12) is easily verified. Thus, let $n < d$.

We may note that by the definition of rectifiability, if $C \subset \mathbb{R}^d$ is closed, then the compact set $S \cap C$ is still countably $\mathcal{H}^n$-rectifiable; besides, (2) holds for all points $x \in S \cap C$ (since it holds for any point $x \in S$). As a consequence, by Theorem 3, we may claim that for any closed subset $C$ of $\mathbb{R}^d$, the following holds:

$$\lim_{r \downarrow 0} \frac{\nu^d((S \cap C)_{\oplus r})}{b_{d-n} r^{d-n}} = \mathcal{H}^n(S \cap C). \tag{13}$$

Let $A$ be as in the assumption.

- Let $\varepsilon > 0$ be fixed. We may observe that the following holds:

$$S_{\oplus r} \cap A \subset (S \cap \operatorname{clos} A)_{\oplus r} \cup (S \cap \operatorname{clos} A_{\oplus \varepsilon} \setminus \operatorname{int} A)_{\oplus r} \qquad \forall r < \varepsilon,$$

where $\operatorname{int} A$ denotes the interior of the set $A$.

Indeed, if $x \in S_{\oplus r} \cap A$, then there exists $y \in S$ with $|x - y| \leq r$ and $y \in \operatorname{clos} A_{\oplus \varepsilon}$. If $x \notin (S \cap \operatorname{clos} A)_{\oplus r}$, we must then have $y \in S \setminus \operatorname{clos} A$, hence $y \in S \cap \operatorname{clos} A_{\oplus \varepsilon} \setminus \operatorname{clos} A$.

By (13), since $\operatorname{clos} A$ and $\operatorname{clos} A_{\oplus \varepsilon} \setminus \operatorname{int} A$ are closed, we have

$$\lim_{r \downarrow 0} \frac{\nu^d((S \cap \operatorname{clos} A)_{\oplus r})}{b_{d-n} r^{d-n}} = \mathcal{H}^n(S \cap \operatorname{clos} A) \stackrel{(11)}{=} \mathcal{H}^n(S \cap A), \tag{14}$$

$$\lim_{r \downarrow 0} \frac{\nu^d(S \cap \operatorname{clos} A_{\oplus \varepsilon} \setminus \operatorname{int} A)_{\oplus r}}{b_{d-n} r^{d-n}} = \mathcal{H}^n(S \cap \operatorname{clos} A_{\oplus \varepsilon} \setminus \operatorname{int} A). \tag{15}$$

Thus,

$$\limsup_{r \downarrow 0} \frac{\nu^d(S_{\oplus r} \cap A)}{b_{d-n} r^{d-n}}$$

$$\leq \limsup_{r\downarrow 0} \frac{\nu^d((S\cap \operatorname{clos} A)_{\oplus r} \cup (S\cap \operatorname{clos} A_{\oplus\varepsilon} \setminus \operatorname{int} A)_{\oplus r})}{b_{d-n}r^{d-n}}$$

$$\leq \limsup_{r\downarrow 0} \frac{\nu^d((S\cap \operatorname{clos} A)_{\oplus r}) + \nu^d((S\cap \operatorname{clos} A_{\oplus\varepsilon} \setminus \operatorname{int} A)_{\oplus r})}{b_{d-n}r^{d-n}}$$

$$\stackrel{(14),(15)}{=} \mathcal{H}^n(S\cap A) + \mathcal{H}^n(S\cap \operatorname{clos} A_{\oplus\varepsilon}\setminus \operatorname{int} A);$$

by taking the limit as $\varepsilon$ tends to 0, we obtain

$$\limsup_{r\downarrow 0} \frac{\nu^d(S_{\oplus r}\cap A)}{b_{d-n}r^{d-n}} \leq \mathcal{H}^n(S\cap A) + \underbrace{\mathcal{H}^n(S\cap \partial A)}_{=0} = \mathcal{H}^n(S\cap A).$$

• Now, let $B$ be a closed set well contained in $A$, that is, the Hausdorff distance between $A$ and $B$ is greater than 0. There then exists $\tilde{r} > 0$ such that $B_{\oplus r} \subset A$, $\forall r < \tilde{r}$. So,

$$\mathcal{H}^n(S\cap B) \stackrel{(13)}{=} \liminf_{r\downarrow 0} \frac{\nu^d((S\cap B)_{\oplus r})}{b_{d-n}r^{d-n}}$$

$$\leq \liminf_{r\downarrow 0} \frac{\nu^d(S_{\oplus r}\cap B_{\oplus r})}{b_{d-n}r^{d-n}}$$

$$\leq \liminf_{r\downarrow 0} \frac{\nu^d(S_{\oplus r}\cap A)}{b_{d-n}r^{d-n}}.$$

Let us consider an increasing sequence of closed sets $\{B_m\}_{m\in\mathbb{N}}$ well contained in $A$ such that $B_m \nearrow \operatorname{int} A$. By taking the limit as $m$ tends to $\infty$, we obtain that

$$\liminf_{r\downarrow 0} \frac{\nu^d(S_{\oplus r}\cap A)}{b_{d-n}r^{d-n}} \geq \lim_{m\to\infty} \mathcal{H}^n(S\cap B_m) = \mathcal{H}^n(S\cap \operatorname{int} A) \stackrel{(11)}{=} \mathcal{H}^n(S\cap A).$$

We summarize:

$$\mathcal{H}^n(S\cap A) \leq \liminf_{r\downarrow 0} \frac{\nu^d(S_{\oplus r}\cap A)}{b_{d-n}r^{d-n}} \leq \limsup_{r\downarrow 0} \frac{\nu^d(S_{\oplus r}\cap A)}{b_{d-n}r^{d-n}} \leq \mathcal{H}^n(S\cap A)$$

and so the thesis follows. □

If we consider the family of random variables $\frac{\nu^d(\Theta_{n_{\oplus r}}\cap A)}{b_{d-n}r^{d-n}}$, for $r$ going to 0, we ask which conditions must be satisfied by a random set $\Theta_n$ so that they are dominated by an integrable random variable. In this way, we could apply the dominated convergence theorem in order to exchange limit and expectation in (12).

**Lemma 7.** *Let $K$ be a compact subset of $\mathbb{R}^d$ and assume that*

$$\eta(B_r(x)) \geq \gamma r^n \qquad \forall x\in K, \forall r\in (0,1) \tag{16}$$

*holds for some $\gamma > 0$ and some probability measure $\eta$ in $\mathbb{R}^d$.*

*Then, for all $r < 2$,*

$$\frac{\nu^d(K_{\oplus r})}{b_{d-n}r^{d-n}} \leq \frac{1}{\gamma}2^n 4^d \frac{b_d}{b_{d-n}}.$$

**Proof.** Since $K_{\oplus r}$ is compact, it is possible to cover it with a finite number $p$ of closed balls $B_{3r}(x_i)$, with $x_i \in K_{\oplus r}$, such that

$$|x_i - x_j| > 3r, \qquad i \neq j.$$

As a consequence, there exist $y_1, \ldots, y_p$ such that:

- $y_i \in K$, $i = 1, \ldots, p$;
- $|y_i - y_j| > r$, $i \neq j$;
- $K_{\oplus r} \subseteq \bigcup_{i=1}^p B_{4r}(y_i)$.

In fact, if $x_i \in K$, then we choose $y_i = x_i$; if $x_i \in K_{\oplus r} \setminus K$, then we choose $y_i \in B_r(x_i) \cap K$. As a consequence, $|y_i - x_i| \leq r$ and $B_{4r}(y_i) \supseteq B_{3r}(x_i)$ for any $i = 1, \ldots, p$. So,

$$\bigcup_{i=1}^p B_{4r}(y_i) \supseteq \bigcup_{i=1}^p B_{3r}(x_i) \supseteq K_{\oplus r}$$

and

$$3r \leq |x_i - x_j| \leq |x_i - y_i| + |y_i - y_j| + |y_j - x_j| \leq 2r + |y_i - y_j|, \qquad i \neq j.$$

For $r < 2$, $B_{r/2}(y_i) \cap B_{r/2}(y_j) = \varnothing$. Since, by hypothesis, $\eta$ is a probability measure satisfying (16), we have that

$$1 \geq \eta\left(\bigcup_{i=1}^p B_{r/2}(y_i)\right) = \sum_{i=1}^p \eta(B_{r/2}(y_i)) \stackrel{(16)}{\geq} p\gamma\left(\frac{r}{2}\right)^n$$

and so

$$p \leq \frac{1}{\gamma}\frac{2^n}{r^n}.$$

In conclusion,

$$\frac{\nu^d(K_{\oplus r})}{b_{d-n}r^{d-n}} \leq \frac{\nu^d(\bigcup_{i=1}^p B_{4r}(y_i))}{b_{d-n}r^{d-n}} \leq \frac{pb_d(4r)^d}{b_{d-n}r^{d-n}} \leq \frac{1}{\gamma}2^n 4^d \frac{b_d}{b_{d-n}}. \qquad \square$$

We are now ready to prove the main result of the paper.

**Proof of Theorem 4.** Let $W$ be a compact window in $\mathbb{R}^d$ and $\Gamma_W(\Theta_n)$ the function defined in (4). Let $Y$ be as in the assumptions of the theorem and $A$ be any Borel subset of $\mathbb{R}^d$ such that

$$A \subset \operatorname{int} W \quad \text{and} \quad \mathbb{E}[\mathcal{H}^n(\Theta_n \cap \partial A)] = 0. \tag{17}$$

Since $\mathbb{E}[Y] < \infty$, we have $Y(\omega) < \infty$ for $\mathbb{P}$-a.e. $\omega \in \Omega$. With no loss of generality, we shall assume that $1/\Gamma_W(\Theta_n) < Y$ a.s.

Let us define

$$\Omega_A := \{\omega \in \Omega : \mathcal{H}^n(\Theta_n(\omega) \cap \partial A) = 0\},$$
$$\Omega_T := \{\omega \in \Omega : \Theta_n(\omega) \text{ is countably } \mathcal{H}^n\text{-rectifiable and closed}\},$$
$$\Omega_Y := \{\omega \in \Omega : Y(\omega) < \infty\},$$
$$\Omega_\Gamma := \left\{\omega \in \Omega : \frac{1}{\Gamma_W(\Theta_n(\omega))} < Y(\omega)\right\};$$

by hypothesis, $\mathbb{P}(\Omega_A) = \mathbb{P}(\Omega_T) = \mathbb{P}(\Omega_Y) = \mathbb{P}(\Omega_\Gamma) = 1$.

Thus, if $\Omega' := \Omega_A \cap \Omega_T \cap \Omega_Y \cap \Omega_\Gamma$, it follows that $\mathbb{P}(\Omega') = 1$.

Let $\omega \in \Omega'$ be fixed. Then:

- $1/Y(\omega) < \Gamma_W(\Theta_n(\omega))$, that is, a probability measure $\eta \ll \mathcal{H}^k$ exists such that

$$\eta(B_r(x)) \geq \frac{1}{Y(\omega)}(\omega) r^n \qquad \forall x \in \Theta_n(\omega) \cap W_{\oplus 1}, \forall r \in (0, 1);$$

- $\mathcal{H}^n(\Theta_n(\omega) \cap \partial A) = 0$.

So, by applying Lemma 6 to $\Theta_n \cap W_{\oplus 1}$ and noting that $(\Theta_n \cap W_{\oplus 1})_{\oplus r} \cap A = \Theta_{n_{\oplus r}} \cap A$, we get

$$\lim_{r \downarrow 0} \frac{\nu^d(\Theta_{n_{\oplus r}}(\omega) \cap A)}{b_{d-n} r^{d-n}} = \mathcal{H}^n(\Theta_n(\omega) \cap A);$$

that is, we may claim that

$$\lim_{r \downarrow 0} \frac{\nu^d(\Theta_{n_{\oplus r}} \cap A)}{b_{d-n} r^{d-n}} = \mathcal{H}^n(\Theta_n \cap A) \qquad \text{almost surely.}$$

Further, for all $\omega \in \Omega'$, $\Theta_n(\omega) \cap W_{\oplus 1}$ satisfies the hypotheses of Lemma 7 and so

$$\frac{\nu^d(\Theta_{n_{\oplus r}}(\omega) \cap A)}{b_{d-n} r^{d-n}} = \frac{\nu^d((\Theta_n(\omega) \cap W_{\oplus 1})_{\oplus r} \cap A)}{b_{d-n} r^{d-n}}$$
$$\leq \frac{\nu^d((\Theta_n(\omega) \cap W_{\oplus 1})_{\oplus r})}{b_{d-n} r^{d-n}} \leq Y(\omega) 2^n 4^d \frac{b_d}{b_{d-n}} \in \mathbb{R}.$$

Let $Z$ be the random variable defined as follows:

$$Z(\omega) := Y(\omega) 2^n 4^d \frac{b_d}{b_{d-n}}, \qquad \omega \in \Omega'.$$

By assumption, $\mathbb{E}[Z] < \infty$, so that the dominated convergence theorem gives Equation (9), and hence the first statement of the theorem.

The second statement is a direct consequence of the first, because if $\Theta_n$ is absolutely continuous in mean, by Definition 1, we have that $\mathbb{E}[\mu_{\Theta_n}](\partial A) = 0$ if $\nu^d(\partial A) = 0$ and

$$\mathbb{E}[\mu_{\Theta_n}](A) = \int_A \lambda_{\Theta_n}(x)\,\mathrm{d}x.$$

Finally, if $\Theta_n$ is stationary, then $\mathbb{P}(x \in \Theta_{n_{\oplus r}})$ is independent of $x$ and the expected measure $\mathbb{E}[\mu_{\Theta_n}]$ is translation invariant, that is, $\lambda_{\Theta_n}(x) = L_{\Theta_n} \in \mathbb{R}_+$ for $\nu^d$-a.e. $x \in \mathbb{R}^d$. It follows that

$$\lim_{r \downarrow 0} \int_A \frac{\mathbb{P}(x \in \Theta_{n_{\oplus r}})}{b_{d-n} r^{d-n}}\,\mathrm{d}x = \lim_{r \downarrow 0} \frac{\mathbb{P}(x_0 \in \Theta_{n_{\oplus r}})}{b_{d-n} r^{d-n}} \nu^d(A)$$

for any $x_0 \in \mathbb{R}^d$ and so (8) follows directly by (7). □

A classical criterion (see, e.g., [11] or [1]) states that $\mu_n$ weakly-* converge to $\mu$ if and only if $\mu_n(A) \to \mu(A)$ for any bounded open set $A$ with $\mu(\partial A) = 0$. Hence, for any random closed set $\Theta_n$ with locally finite $n$-dimensional measure a.s. satisfying the assumptions of Theorem 4, it follows that the measures $\mu^{\oplus r}$ defined by

$$\mu^{\oplus r}(A) := \int_A \frac{\mathbb{P}(x \in \Theta_{n_{\oplus r}})}{b_{d-n} r^{d-n}}\,\mathrm{d}x \qquad \forall A \in \mathcal{B}_{\mathbb{R}^d},$$

weakly-* converge to the measure $\mathbb{E}[\mu_{\Theta_n}]$ as $r \downarrow 0$. As a consequence, the associated density functions

$$\delta_n^{\oplus r}(x) := \frac{\mathbb{P}(x \in \Theta_{n_{\oplus r}})}{b_{d-n} r^{d-n}}$$

can be taken to constitute a *weak approximating family* of the density $\delta \mathbb{E}[\mu_{\Theta_n}]$ of $\Theta_n$ when the process is absolutely continuous in mean.

*Remark 8 (Mean density as a pointwise limit).* It is tempting to try to exchange limit and integral in (7), to obtain

$$\lim_{r \downarrow 0} \frac{\mathbb{P}(x \in \Theta_{n_{\oplus r}})}{b_{d-n} r^{d-n}} = \lambda_{\Theta_n}(x), \tag{18}$$

at least for $\nu^d$-a.e. $x \in \mathbb{R}^d$. The proof of the validity of this formula for absolutely continuous (in mean) processes seems to be quite a delicate problem, with the only exception being stationary processes; this is the main reason why the existing literature has extensively considered this case [3]; we nonetheless wish to stress that global stationarity is only a sufficient condition for (18) to hold. Indeed, we already know that in the extreme cases $n = d$ and $n = 0$, it is not hard to prove it.

In the case $n = d$, we know from Remark 2 that $\lambda_{\Theta_d}(x) = \mathbb{P}(x \in \Theta_d)$ for $\nu^d$-a.e. $x$, and obviously $\mathbb{P}(x \in \Theta_{d_{\oplus r}})$ converges to $\mathbb{P}(x \in \Theta_d)$ for all $x$.

In the case $n = 0$, let $\Theta_0 = \{X\}$, with $X$ an absolutely continuous random point in $\mathbb{R}^d$ with pdf $f_X$, and note that

$$\frac{\mathbb{P}(x \in \Theta_{0 \oplus r})}{b_d r^d} = \frac{\mathbb{P}(X \in B_r(x))}{b_d r^d} = \frac{1}{b_d r^d} \int_{B_r(x)} f_X(y) \, dy.$$

Therefore, (18) with $\lambda_{\Theta_0} = f_X$ holds at any Lebesgue point of $f_X$ (i.e., a point $x$ such that the mean values of $y \mapsto |f_X(y) - f_X(x)|$ on $B_r(x)$ are infinitesimal as $r \downarrow 0$) and hence for $\nu^d$-a.e. $x \in \mathbb{R}^d$. This is the main reason why the problem of local density estimation has already been solved [10, 22].

We may note that

$$\mathbb{P}(x \in \Theta_{n \oplus r}) = \mathbb{P}(\Theta_n \cap B_r(x) \neq \varnothing) = T_{\Theta_n}(B_r(x)),$$

thus making explicit the connection to $T_{\Theta_n}$, the capacity functional characterizing the probability law of the random set $\Theta_n$ [16]. Therefore, the family of functions $\delta^{\oplus r}$ may suggest estimators of the mean density of $\Theta_n$ in terms of the empirical capacity functional of $\Theta_n$.

For instance, even if the special case $n = 0$ can be handled with much more elementary tools, we may note how the well-known and well-studied case of a random point $X$ in $\mathbb{R}^d$ is consistent with our framework. It is immediate to check that $X$ satisfies the hypotheses of Theorem 4 (it is sufficient to choose $\eta := \mathcal{H}^0(X(\omega) \cap \cdot)$) and that, in the particular case $d = 1$ with $X$ absolutely continuous r.v. with pdf $f_X$, (18) becomes

$$\lim_{r \downarrow 0} \frac{\mathbb{P}(X \in [x - r, x + r])}{2r} = f_X,$$

which leads to the usual probability density estimation by histograms.

## 4. Applications

In many real applications, $\Theta_n$ is given by a random collection of geometrical objects, so that it may be described as the union of a family of random closed sets $E_i$ in $\mathbb{R}^d$,

$$\Theta_n = \bigcup_i E_i, \tag{19}$$

and a problem of interest is to determine or estimate the mean density $\lambda_{\Theta_n}$ of $\Theta_n$. Such unions of random sets can be taken as models for the so-called *particle processes*, widely studied in stochastic geometry.

We may note that Theorem 4 seems to require sufficient regularity of the $E_i$'s, rather than specific assumptions about their probability law or the stochastic dependence among them. In the present section, we provide concrete examples and results in this direction, in order to clarify the wide range of applicability of our main theorem and to show how some known results in the current literature for stationary random closed sets follow from (8).

For the sake of simplicity, we start by considering the case in which the union (19) of the $E_i$ is finite a.s. so that there exists a positive integer-valued random variable $\Phi$ such that

$$\Theta_n = \bigcup_{i=1}^{\Phi} E_i. \tag{20}$$

In the next example, we consider a class of random sets of this kind, which might be taken to model a class of *birth-and-growth stochastic processes*, relevant in real applications, as discussed in Example 2.

*Example 1.* A class of random sets satisfying the hypotheses of Theorem 4 is given by all sets $\Theta_n$ which are a random union of countably $\mathcal{H}^n$-rectifiable random closed sets in $\mathbb{R}^d$, as in (20), such that

(i) $\mathbb{E}[\Phi] < \infty$;
(ii) $E_1, E_2, \ldots$ are i.i.d. as $E$ and independent of $\Phi$;
(iii) $\mathbb{E}[\mathcal{H}^n(E)] = C < \infty$ and $\exists \gamma > 0$ such that for any $\omega \in \Omega$,

$$\mathcal{H}^n(E(\omega) \cap B_r(x)) \geq \gamma r^n \qquad \forall x \in E(\omega), \forall r \in (0,1). \tag{21}$$

We can choose $\eta(\cdot) := \frac{\mathcal{H}^n(\Theta_n(\omega) \cap \cdot)}{\mathcal{H}^n(\Theta_n(\omega))}$ for any fixed $\omega \in \Omega$. As a consequence, $\eta$ is a probability measure absolutely continuous with respect to $\mathcal{H}^n$ and such that

$$\eta(B_r(x)) \geq \frac{\gamma}{\mathcal{H}^n(\Theta(\omega))} r^n \qquad \forall x \in \Theta_n(\omega),\ r \in (0,1).$$

In fact, if $x \in \Theta_n(\omega)$, then there exists an $\bar{\imath}$ such that $x \in E_{\bar{\imath}}(\omega)$; since $\Theta_n(\omega) = \bigcup_{i=1}^{\Phi(\omega)} E_i(\omega)$, we have

$$\eta(B_r(x)) = \frac{\mathcal{H}^n(\Theta_n(\omega) \cap B_r(x))}{\mathcal{H}^n(\Theta_n(\omega))} \geq \frac{\mathcal{H}^n(E_{\bar{\imath}}(\omega) \cap B_r(x))}{\mathcal{H}^n(\Theta_n(\omega))} \geq \frac{\gamma}{\mathcal{H}^n(\Theta(\omega))} r^n.$$

As a result, the function $\Gamma(\Theta_n)$ defined in (4) is such that $\Gamma(\Theta_n) \geq \gamma/\mathcal{H}^n(\Theta_n) =: 1/Y$ and so it remains to verify only that $\mathbb{E}[\mathcal{H}^n(\Theta_n)] < \infty$:

$$\begin{aligned}
\mathbb{E}[\mathcal{H}^n(\Theta_n)] &= \mathbb{E}[\mathbb{E}[\mathcal{H}^n(\Theta_n)|\Phi]] \\
&= \sum_{k=1}^{\infty} \mathbb{E}\left[\mathcal{H}^n\left(\bigcup_{i=1}^{k} E_i\right)\bigg| \Phi = k\right] \mathbb{P}(\Phi = k) \\
&\stackrel{\text{(ii)}}{\leq} \sum_{k=1}^{\infty} \sum_{i=1}^{k} \mathbb{E}[\mathcal{H}^n(E_i)] \mathbb{P}(\Phi = k)
\end{aligned}$$

$$\stackrel{\text{(iii)}}{=} \sum_{k=1}^{\infty} Ck\mathbb{P}(\Phi = k)$$

$$= C\mathbb{E}[\Phi] \stackrel{\text{(i)}}{<} \infty.$$

Note that we have not made any particular assumption on the probability laws of $\Phi$ and $E$. Further, it is clear that the same proof holds, even in the case in which the $E_i$'s are not i.i.d., provided that $\mathbb{E}[\mathcal{H}^n(E_i)] \leq C$ for all $i$ and (21) is true for any $E_i$ (with $\gamma$ independent of $\omega$ and $i$).

In the next example, we show how unions of time-dependent random closed sets may be used to model a class of time-dependent geometric processes as well; it also provides an example of an absolutely continuous in mean, but not stationary, random closed set.

*Example 2 (Birth-and-growth process).* Roughly speaking, a *birth-and-growth process* is a family $\{\Theta^t\}_t$ of random closed sets which develop in time, according to a given growth model, from points (*nuclei*, or *germs*) that are born at random both in space and time (for details, see [6] and references therein). These kinds of processes are described by dynamic germ-grain models (see, e.g., [9, 24]), that is, once born, each germ generates a grain subject to surface growth with a speed $G$ which is, in general, space-time dependent. The nucleation process $\{T_n, X_n\}_{n \in \mathbb{N}}$, where $T_n$ is the $\mathbb{R}_+$-valued random variable representing the time of birth of the $n$th nucleus and $X_n$ is the $\mathbb{R}^d$-valued random variable representing the spatial location of the nucleus born at time $T_n$, is commonly modelled by a marked point process $N$. In particular, $N$ is a random counting measure on $\mathbb{R}_+ \times \mathbb{R}^d$ such that, for any fixed $t \in \mathbb{R}_+$, the random number $N([0,t] \times \mathbb{R}^d)$ of nuclei born up to time $t$ is finite with probability 1.

For the sake of simplicity, we assume here that $N$ is given by an inhomogeneous Poisson point process in $\mathbb{R}_+ \times \mathbb{R}^d$ with intensity $\alpha(t,x)$ (i.e., $\mathbb{E}[N(\mathrm{d}(t,x)] = \alpha(t,x)\,\mathrm{d}t\,\mathrm{d}x$) and that the growth occurs with a constant normal velocity $G > 0$ so that, for any fixed time $t$, $\Theta^t$ is a finite union of random balls in $\mathbb{R}^d$:

$$\Theta^t = \bigcup_{i: T_i \leq t} B_{G(t-T_i)}(X_i).$$

Considering now the random closed set $\partial \Theta^t$, we observe that it is absolutely continuous in mean, but not stationary, for any $t > 0$. We may prove this (by contradiction) as follows.

Assume $\mathbb{E}[\mathcal{H}^{d-1}(\partial \Theta^t \cap \cdot)]$ to be not absolutely continuous with respect to $\nu^d$; there then exists some $A \subset \mathbb{R}^d$ with $\nu^d(A) = 0$ such that $\mathbb{E}[\mathcal{H}^{d-1}(\partial \Theta^t \cap A)] > 0$.

It is clear that

$$\mathbb{E}[\mathcal{H}^{d-1}(\partial \Theta^t \cap A)] > 0 \quad \Rightarrow \quad \mathbb{P}(\mathcal{H}^{d-1}(\partial \Theta^t \cap A) > 0) > 0$$

and

$$\mathbb{P}(\mathcal{H}^{d-1}(\partial \Theta^t \cap A) > 0) \leq \mathbb{P}(\exists (T_j, X_j) : \mathcal{H}^{d-1}(\partial B_{G(t-T_j)}(X_j) \cap A) > 0).$$

As a consequence, we have

$$\mathbb{E}[\mathcal{H}^{d-1}(\partial \Theta^t \cap A)] > 0 \quad \Rightarrow \quad \mathbb{P}(\Phi(\mathcal{A}) \neq 0) > 0,$$

where

$$\mathcal{A} := \{(s,y) \in [0,t] \times \mathbb{R}^d : \mathcal{H}^{d-1}(\partial B_{G(t-s)}(y) \cap A) > 0\}.$$

Denoting by $\mathcal{A}_s := \{y \in \mathbb{R}^d : (s,y) \in \mathcal{A}\}$ the section of $\mathcal{A}$ at time $s$ and by $\mathcal{A}^y := \{s > 0 : (s,y) \in \mathcal{A}\}$ the section of $\mathcal{A}$ at $y$, we note that $\nu^1(\mathcal{A}^y) = 0$ for all $y$ because $\nu^d(A) = 0$ (it suffices to use spherical coordinates centered at $y$ to obtain that the $\nu^1$-a.e. ball with radius $s$ centered at $y$ intersects $A$ in an $\mathcal{H}^{d-1}$-negligible set). Therefore, we may apply Fubini's theorem to get

$$\int_0^\infty \nu^d(\mathcal{A}_s) \, \mathrm{d}s = \int_0^\infty \int_{\mathbb{R}^d} \chi_{\mathcal{A}_s} \, \mathrm{d}y \, \mathrm{d}s = \int_{\mathbb{R}^d} \int_0^\infty \chi_{\mathcal{A}^y} \, \mathrm{d}s \, \mathrm{d}y = \int_{\mathbb{R}^d} \nu^1(\mathcal{A}^y) \, \mathrm{d}y = 0.$$

It follows that $\nu^d(\mathcal{A}_s) = 0$ for $\nu^1$-a.s. $s \in [0,t]$ and so

$$\mathbb{E}[\Phi(\mathcal{A})] = \int_\mathcal{A} \alpha(s,y) \, \mathrm{d}s \, \mathrm{d}y = \int_0^t \int_{\mathcal{A}_s} \alpha(s,y) \, \mathrm{d}y \, \mathrm{d}s = 0.$$

But this is an absurd, since

$$\mathbb{P}(\Phi(\mathcal{A}) \neq 0) > 0 \Rightarrow \mathbb{E}[\Phi(\mathcal{A})] > 0.$$

Observing that Theorem 4 applies to $\partial \Theta^t$ (see [25], Chapter 4 for details), as a consequence of (7), we have

$$\lim_{r \downarrow 0} \int_A \frac{\mathbb{P}(x \in \partial \Theta^t_{\oplus r})}{2r} \, \mathrm{d}x = \int_A \lambda_{\partial \Theta^t}(x) \, \mathrm{d}x$$

for any bounded Borel set $A \subset \mathbb{R}^d$ with $\nu^d(\partial A) = 0$, where $\lambda_{\partial \Theta^t}$ is also known as the *mean surface density* associated with the birth-and-growth process $\{\Theta^t\}_t$. Mean volume and surface densities of birth-and-growth processes are of a great interest in applications (see [6] and references therein).

Coming back to consider the basic model described in Example 1, we now prove a result which relates the probability that a point $x$ belongs to the set $\Theta_{n_{\oplus r}}$ to the mean number of $E_i$ which intersect the ball $B_r(x)$, and which might therefore be useful in statistical applications.

**Proposition 9.** *Let $n < d$, let $\Phi$ be a positive integer-valued random variable with $\mathbb{E}[\Phi] < \infty$ and let $\{E_i\}$ be a collection of random closed sets with locally finite $n$-dimensional measure a.s. Let $\Theta_n$ be the random closed set defined as*

$$\Theta_n = \bigcup_{i=1}^\Phi E_i.$$

*Approximation of mean densities*  1237

If $E_1, E_2, \ldots$ are i.i.d. as $E$ and independent of $\Phi$, then, for any $x \in \mathbb{R}^d$ such that $\mathbb{P}(x \in E) = 0$,

$$\lim_{r\downarrow 0} \frac{\mathbb{P}(x \in \Theta_{n_{\oplus r}})}{b_{d-n}r^{d-n}} = \lim_{r\downarrow 0} \frac{\mathbb{E}[\#\{E_i : x \in E_{i_{\oplus r}}\}]}{b_{d-n}r^{d-n}} = \mathbb{E}[\Phi]\lim_{r\downarrow 0} \frac{\mathbb{P}(x \in E_{\oplus r})}{b_{d-n}r^{d-n}}, \tag{22}$$

*provided that at least one of the above limits exists.*

**Proof.** The following chain of equalities holds:

$$\begin{aligned}
\mathbb{P}(x \in \Theta_{n_{\oplus r}}) &= \mathbb{P}\left(x \in \bigcup_{i=1}^{\Phi} E_{i_{\oplus r}}\right) \\
&= 1 - \mathbb{P}\left(x \notin \bigcup_{i=1}^{\Phi} E_{i_{\oplus r}}\right) \\
&= 1 - \mathbb{P}\left(\bigcap_{i=1}^{\Phi} \{x \notin E_{i_{\oplus r}}\}\right) \\
&= 1 - \sum_{k=1}^{\infty} \mathbb{P}\left(\bigcap_{i=1}^{k} \{x \notin E_{i_{\oplus r}}\}\Big|\Phi = k\right)\mathbb{P}(\Phi = k) \\
&= 1 - \sum_{k=1}^{\infty} [\mathbb{P}(x \notin E_{\oplus r})]^k \mathbb{P}(\Phi = k),
\end{aligned}$$

the $E_i$'s being i.i.d. and independent of $\Phi$. Denoting by G the probability generating function of the random variable $\Phi$ and $s(r) := \mathbb{P}(x \in E_{\oplus r})$, we get that

$$\mathbb{P}(x \in \Theta_{n_{\oplus r}}) = 1 - \mathrm{G}(1 - s(r)). \tag{23}$$

We now observe that

$$\begin{aligned}
\mathbb{E}[\#\{E_i : x \in E_{i_{\oplus r}}\}] &= \sum_{k=1}^{\infty} \mathbb{E}\left[\sum_{i=0}^{k} \mathbf{1}_{E_{i_{\oplus r}}}(x)\Big|\Phi = k\right]\mathbb{P}(\Phi = k) \\
&= \sum_{k=1}^{\infty} k\mathbb{P}(x \in E_{\oplus r})\mathbb{P}(\Phi = k) \tag{24} \\
&= \mathbb{E}[\Phi]\mathbb{P}(x \in E_{\oplus r})
\end{aligned}$$

and recall that $\mathbb{E}[\Phi] = \mathrm{G}'(1)$ and $1 = \mathrm{G}(1)$. By hypothesis, we know that $s(r) \downarrow 0$ as $r \downarrow 0$. Thus, by (23) and (24), we have that if $s(r) = 0$ for some $r > 0$, then $\mathbb{P}(x \in \Theta_{n_{\oplus R}}) = \mathbb{E}[\#\{E_i : x \in E_{i_{\oplus R}}\}] = 0$ for all $\leq R$ and so the assertion trivially follows, whereas if

$s(r) > 0$ for all $r > 0$, we have

$$\lim_{r \downarrow 0} \frac{\mathbb{P}(x \in \Theta_{n_{\oplus r}})}{\mathbb{E}[\#\{E_i : x \in E_{i_{\oplus r}}\}]} = \lim_{r \downarrow 0} \frac{G(1) - G(1 - s(r))}{G'(1) s(r)}$$
$$= \frac{1}{G'(1)} \lim_{r \downarrow 0} \frac{G(1 - s(r)) - G(1)}{-s(r)} \quad (25)$$
$$= \frac{1}{G'(1)} G'(1) = 1.$$

In conclusion, we obtain

$$\lim_{r \downarrow 0} \frac{\mathbb{P}(x \in \Theta_{n_{\oplus r}})}{b_{d-n} r^{d-n}} = \lim_{r \downarrow 0} \frac{\mathbb{P}(x \in \Theta_{n_{\oplus r}})}{\mathbb{E}[\#\{E_i : x \in E_{i_{\oplus r}}\}]} \frac{\mathbb{E}[\#\{E_i : x \in E_{i_{\oplus r}}\}]}{b_{d-n} r^{d-n}}$$
$$\stackrel{(25)}{=} \lim_{r \downarrow 0} \frac{\mathbb{E}[\#\{E_i : x \in E_{i_{\oplus r}}\}]}{b_{d-n} r^{d-n}} \quad (26)$$
$$\stackrel{(24)}{=} \mathbb{E}[\Phi] \lim_{r \downarrow 0} \frac{\mathbb{P}(x \in E_{\oplus r})}{b_{d-n} r^{d-n}}. \quad (27)$$

□

Note that for any random closed set $\Theta_n$ as in Proposition (9), by (22), we infer that the probability that a point $x$ belongs to the intersection of two or more enlarged sets $E_i$ is infinitesimally faster than $r^{d-n}$. In fact, denoting by $F_r(x)$ this event and by $\chi_A : \Omega \to \{0, 1\}$ the characteristic function of an event $A$, we have

$$\chi_{F_r(x)} \leq \sum_i \chi_{\{x \in E_{i \oplus r}\}} - \chi_{\{x \in \Theta_{n \oplus r}\}},$$

so that, since (24) gives

$$\mathbb{E}\left[\sum_i \chi_{\{x \in E_{i \oplus r}\}}\right] = \mathbb{E}[\Phi] \mathbb{P}(x \in E_{\oplus r}),$$

by taking expectations on both sides and dividing by $r^{d-n}$, we get $\mathbb{P}(F_r(x))/r^{d-n} \to 0$ as $r \downarrow 0$.

**Remark 10.** Whenever $\Theta_n$ is absolutely continuous in mean and it is possible to exchange limit and integral in (7), we can use the fact that $A$ is arbitrary to obtain

$$\lambda_{\Theta_n}(x) = \mathbb{E}[\Phi] \lim_{r \downarrow 0} \frac{\mathbb{P}(x \in E_{\oplus r})}{b_{d-n} r^{d-n}} = \mathbb{E}[\Phi] \lambda_E(x)$$

for $\nu^d$-a.e. $x \in \mathbb{R}^d$, where $\lambda_{\Theta_n}$ and $\lambda_E$ are the densities of $\mathbb{E}[\mu_{\Theta_n}]$ and $\mathbb{E}[\mu_E]$, respectively. In particular, when $E$ is stationary (which implies that $\Theta_n$ is also stationary), $\lambda_{\Theta_n}(x) \equiv$

$L_{\Theta_n} \in \mathbb{R}^+$ and $\lambda_E(x) \equiv L_E \in \mathbb{R}^+$ so that

$$L_{\Theta_n} = \mathbb{E}[\Phi] \lim_{r \downarrow 0} \frac{\mathbb{P}(x_0 \in E_{\oplus r})}{b_{d-n} r^{d-n}} = \mathbb{E}[\Phi] L_E.$$

Note that, in general, $\Theta_n(\omega)$ may be unbounded because the $E_i(\omega)$'s are not compact (e.g., $E_i$ may be a random line) or the union in not finite. Random closed sets of this type are often taken to model real problems in applications so that in the literature, many geometric processes like these are investigated, as point-, line-, segment- or plane processes, random mosaics, grain processes, etcetera (see, e.g., [3, 24]). In any case, it is commonly assumed in stochastic geometry that $\Theta_n$ is locally finite (i.e., the number of $E_i$ hitting any compact subset of $\mathbb{R}^d$ is finite a.s.); this is the reason why the compact window $W \subset \mathbb{R}^d$ is introduced in Theorem 4. Indeed, denoting by $\Phi^W$ the random number of $E_i$'s hitting $W$ and considering the restriction of $\Theta_n$ to $W$, we may represent it as a finite union of random compact sets, almost surely, (and thus as in (20)) in the following way

$$\Theta_n \cap W = \bigcup_{i=1}^{\Phi^W} E_i^W,$$

where $E_i^W = E_i \cap W$.

In the case of unbounded random sets, we consider the so-called *Poisson line process*, as a simple example of applicability of Theorem 4.

***Example 3 (Poisson line process).*** Let $\Theta_1$ be the random closed set of Hausdorff dimension 1 a.s. associated with a Poisson line process in the plane (see, e.g., [24] for details). $\Theta_1$ is then a random collection of lines in $\mathbb{R}^2$ such that the mean number of lines hitting any compact planar set is finite. Without any further requirement on $\Theta_1$ (e.g., stationarity or isotropy), it is easy to see that the hypotheses of Theorem 4 are satisfied with $Y = \mathcal{H}^1(\Theta_1 \cap W_{\oplus 1})$. Indeed, by choosing

$$\eta(\cdot) := \frac{\mathcal{H}^1(\Theta_1(\omega) \cap W_{\oplus 1} \cap \cdot)}{\mathcal{H}^1(\Theta_1(\omega) \cap W_{\oplus 1})} \qquad \forall \omega \in \Omega, \tag{28}$$

it follows that

$$\eta(B_r(x)) \geq \frac{1}{\mathcal{H}^1(\Theta_1(\omega) \cap W_{\oplus 1})} r \qquad \forall x \in \Theta_1(\omega) \cap W_{\oplus 1}, \forall r \in (0,1),$$

and, denoting by $\Phi^W$ the random number of lines hitting $W_{\oplus 1}$,

$$\mathbb{E}[\mathcal{H}^1(\Theta_1 \cap W_{\oplus 1})] \leq C \mathbb{E}[\Phi^W] < \infty,$$

for a suitable constant $C$ which depends on the dimension of the window $W$.

Well known in the literature is the particular case in which $\Theta_1$ is stationary and isotropic. In this case, $\Theta_1$ can be described by a point process in the cylinder $\mathbf{C}^* =$ $\{(\cos\alpha, \sin\alpha,\ p) : p \in \mathbb{R}, \alpha \in (0,\pi]\}$ in $\mathbb{R}^3$, with intensity measure $\Lambda(\mathrm{d}(p,\alpha)) = L \cdot \mathrm{d}p \cdot \frac{\mathrm{d}\alpha}{2\pi}$ for some constant $L > 0$. It is proved that the mean density $\lambda_{\Theta_1}$ of $\Theta_1$ is equal to $L$ and that the number $N_r$ of lines of $\Theta_1$ hitting the ball $B_r(0)$ is a Poisson random variable with mean $2rL$ (see [24], pages 249–250). We may note that the same result for the mean density can also be obtained as an application of Theorem 4: by (8), we have that

$$\lambda_{\Theta_1} = \lim_{r\downarrow 0} \frac{\mathbb{P}(\Xi \cap B_r(0) \neq \varnothing)}{2r} = \lim_{r\downarrow 0} \frac{\mathbb{P}(N_r \geq 1)}{2r} = \lim_{r\downarrow 0} \frac{1 - \mathrm{e}^{-2rL_\Xi}}{2r} = L.$$

Let us observe that in the above example, in order to apply Theorem 4, for all $\omega \in \Omega$, we have chosen the normalized measure $\eta$ as the restriction of the Hausdorff measure $\mathcal{H}^1$ to $\Theta_1(\omega) \cap W_{\oplus 1}$. Indeed, problems may arise in identifying a measure $\eta$ needed for the application of the quoted theorem; in some cases, a proper choice of $\eta$ can be made by referring to another suitable random set which contains the relevant random set $\Theta_n$. We further clarify this procedure by means of the following example, which suggests, in particular, a method to study more complex random sets, as briefly discussed at the end of the section.

***Example 4 (Segment process).*** Let $\Theta_1 = \bigcup_i S_i$ be a collection of i.i.d. random segments $S_i$ in $\mathbb{R}^d$ with finite expected length, such that the mean number of segments hitting any compact set is finite.

Given a compact window $W \subset \mathbb{R}^d$, if, for all $\omega \in \Omega$, we choose $\eta$ as in (28), then it is easy to check that $\Gamma_W(\Theta_1) \geq \min\{1, L\}/\mathcal{H}^1(\Theta_1 \cap W_{\oplus 1}) =: 1/Y$, where

$$L := \min_{i : S_i \cap W_{\oplus 1} \neq \varnothing} \{\mathcal{H}^1(S_i)\}.$$

With respect to Theorem 4, the above is not a good choice for $\eta$, because, unless $L$ is bounded from below, we may well have $\mathbb{E}[Y] = \infty$; in this case, a possible solution to the problem is to extend all the segments having length less than 2 (the extension can be done homothetically from the center of the segment, so that measurability of the process is preserved). In particular, for any $\omega \in \Omega$, let

$$\widetilde{S}_i(\omega) = \begin{cases} S_i(\omega), & \text{if } \mathcal{H}^1(S_i(\omega)) \geq 2, \\ S_i(\omega) \text{ extended to length 2}, & \text{if } \mathcal{H}^1(S_i(\omega)) < 2, \end{cases}$$

and $\widetilde{\Theta}_1 := \bigcup_i \widetilde{S}_i$. In this way, Theorem 4 applies with $Y := \mathcal{H}^1(\widetilde{\Theta}_1 \cap W_{\oplus 1})$, by now choosing $\eta$ as in (28) and replacing $\Theta_1 \cap W_{\oplus 1}$ with $\widetilde{\Theta}_1 \cap W_{\oplus 1}$ (see [25], Chapter 4 for details).

***Remark 11.*** A well-studied segment process in the literature is the so-called *stationary Poisson segment process* in $\mathbb{R}^d$ (see [3, 24]). By proceeding as in the above example, the known result concerning the mean density of this kind of process can be reobtained by applying (8) in Theorem 4.

Example 4 provides a method of applicability of Theorem 4 to more complex random closed sets $\Theta_n$ with Hausdorff dimension $1 \leq n < d$ a.s. (see also [25], Chapter 4 for an application to Boolean models of spheres). An important class of 1-dimensional random sets are the so-called *fiber processes* (i.e., random collections of 1-rectifiable curves [3]). We recall that a set $S \subset \mathbb{R}^d$ is said to be $n$-rectifiable if it is representable as the image of a compact set $K \subset \mathbb{R}^n$, with $f \colon \mathbb{R}^n \to \mathbb{R}^d$ Lipschitz, and we point out that condition 2 is satisfied with $\eta(\cdot) = \mathcal{H}^n(\widetilde{S} \cap \cdot)$ for some closed set $\widetilde{S} \supseteq S$ if $f$ admits a Lipschitz inverse (see [1], page 111). Hence, if $\Theta_1$ is a sufficiently regular fiber process, an argument similar to that in the above example might be applied, by considering as $\widetilde{\Theta}_1$ the random closed set given by the union of suitably extended fibers, so that information about the measure $\mathbb{E}[\mu_{\Theta_1}]$ might also be obtained under hypotheses of inhomogeneity of the process. A relevant real system which can be modelled as a fiber process is the system of vessels in tumor-driven angiogenesis; estimation of the mean length intensity of such a system is useful for suggesting important methods of diagnosis and of dose response in clinical treatments [5, 8].

## Acknowledgements

It is a pleasure to acknowledge fruitful discussions with A. Micheletti of Milan University. We also thank the referees for their detailed and helpful comments. VC acknowledges the warm hospitality of the Austrian Academy of Sciences at RICAM (Radon Institute for Computational and Applied Mathematics) in Linz, chaired by Professor H. Engl.

## References

[1] Ambrosio, L., Fusco, N. and Pallara, D. (2000). *Functions of Bounded Variation and Free Discontinuity Problems*. Oxford: Clarendon Press. MR1857292

[2] Baddeley, A.J. and Molchanov, I.S. (1997). On the expected measure of a random set. In *Proceedings of the International Symposium on Advances in Theory and Applications of Random Sets (Fontainebleau, 1996)* 3–20. River Edge, NJ: World Scientific. MR1654394

[3] Beneš, V. and Rataj, J. (2004). *Stochastic Geometry: Selected Topics*. Boston: Kluwer. MR2068607

[4] Capasso, V. (ed.) (2003). *Mathematical Modelling for Polymer Processing. Polymerization, Crystallization, Manufacturing. Springer Series on Mathematics in Industry* **2**. Heidelberg: Springer-Verlag. MR1964653

[5] Capasso, V. and Micheletti, A. (2005). Stochastic geometry and related statistical problems in biomedicine. In *Complex Systems in Biomedicine* (A. Quarteroni et al., eds.). Milano: Springer.

[6] Capasso, V. and Villa, E. (2007). On mean densities of inhomogeneous geometric processes arising in material science and medicine. *Image Anal. Stereol.* **26** 23–36. MR2337324

[7] Capasso, V. and Villa, E. (2008). On the geometric densities of random closed sets. *Stoch. Anal. Appl.* **26** 784–808. MR2431079

[8] Carmeliet, P. and Jain, R.K. (2000). Angiogenesis in cancer and other diseases. *Nature* **407** 249–257.

[9] Daley, D.J. and Vere-Jones, D. (1998). *An Introduction to the Theory of Point Processes*. New York: Springer. MR0950166

[10] Devroye, L. (1987). *A Course in Density Estimation*. Boston: Birkhauser Verlag. MR0891874

[11] Evans, L.C. and Gariepy, R.F. (1992). *Measure Theory and Fine Properties of Functions*. Boca Raton: CRC Press. MR1158660

[12] Hahn, U., Micheletti, A., Pohlink, R., Stoyan, D. and Wendrock, H. (1999). Stereological analysis and modeling of gradient structures. *J. Microsc.* **195** 113–124.

[13] Jeulin, D. (2002). Modelling random media. *Image Anal. Stereol.* **21** (Suppl. 1) S31–S40.

[14] Kärkkäinen, S., Jensen, E.B.V. and Jeulin, D. (2002). On the orientational analysis of planar fibre systems from digital images. *J. Microsc.* **207** 69–77. MR1913181

[15] Kolmogorov, A.N. (1956). *Foundations of the Theory of Probability*, 2nd English ed. New York: Chelsea. MR0079843

[16] Matheron, G. (1975). *Random Sets and Integral Geometry*. New York: Wiley. MR0385969

[17] Miles, R.E. and Serra, J. (eds.) (1978). *Geometrical Probability and Biological Structures: Buffon's 200th Anniversary. Lecture Notes in Biomathematics* **23**. Berlin–New York: Springer-Verlag. MR0518158

[18] Møller, J. (1992). Random Johnson–Mehl tessellations. *Adv. in Appl. Probab.* **24** 814–844. MR1188954

[19] Møller, J. (1994). *Lectures on Random Voronoi Tessellations. Lecture Notes in Statistics* **87**. New York–Berlin–Heidelberg: Springer-Verlag. MR1295245

[20] Pestman, W.R. (1998). *Mathematical Statistics: An Introduction*. Berlin: de Gruyter. MR1643593

[21] Robbins, H.E. (1944). On the measure of a random set. *Ann. Math. Statist.* **15** 70–74. MR0010347

[22] Silverman, B.W. (1986). *Density Estimation for Statistics and Data Analysis*. London: Chapman & Hall. MR0848134

[23] Serra, J. (1984). *Image Analysis and Mathematical Morphology*. London: Academic Press. MR0753649

[24] Stoyan, D., Kendall, W.S. and Mecke, J. (1995). *Stochastic Geometry and Its Applications*. Chichester: Wiley. MR0895588

[25] Villa, E. (2007). Methods of geometric measure theory in stochastic geometry. Ph.D. thesis, University of Milan, Milan.

[26] Zähle, M. (1982). Random processes of Hausdorff rectifiable closed sets. *Math. Nachr.* **108** 49–72. MR0695116